\documentstyle[12pt,bezier]{article}
\newtheorem{theorem}{Theorem}
\newtheorem{definition}[theorem]{Definition}
\newtheorem{corollary}[theorem]{Corollary}
\newtheorem{conjecture}[theorem]{Conjecture}
\newtheorem{example}[theorem]{Example}
\newtheorem{lemma}[theorem]{Lemma}

\newtheorem{proposition}[theorem]{Proposition}
\newtheorem{remark}[theorem]{Remark}
\newtheorem{problem}[theorem]{Problem}

\font\tengoth=eufm10 scaled \magstep1
\font\sevengoth=eufm7 scaled \magstep1
\font\fivegoth=eufm5 scaled \magstep1
\newfam\gothfam \def\frak{\fam\gothfam\tengoth}
\textfont\gothfam=\tengoth
\scriptfont\gothfam=\sevengoth
\scriptscriptfont\gothfam=\fivegoth

\def\boldface{\sf}

\catcode`\@=11

\def\ps@myheadings{\let\@mkboth\@gobbletwo
\def\@oddhead{\ifnum\count0=1 \hfill\else
\rightmark \hfil \rm\thepage\fi}%
\def\@oddfoot{\ifnum\count0=1 \hfill \rm 1 \hfill \else
\hfill\fi}
\def\@evenhead%
{\rm\leftmark\hfil\rm\thepage}%
\def\@evenfoot{}\def\sectionmark##1{}
\def\subsectionmark##1{}}

\def\@begintheorem#1#2{\it \trivlist \item[\hskip
 \labelsep{\boldface #1\ #2.}]}
\def\@opargbegintheorem#1#2#3{\it \trivlist\item[\hskip%
 \labelsep{\boldface #1\ #2.\ (#3)}]}
\def\@endtheorem{\endtrivlist}

\def\@listI{\leftmargin\leftmargini \parsep 1pt plus 2.5pt
 minus 1pt\topsep 5pt plus 2pt minus 3pt%
 \itemsep 0pt plus 2.5pt minus 1pt}
\let\@listi\@listI
\@listi

\def\@sect#1#2#3#4#5#6[#7]#8{\ifnum #2>\c@secnumdepth%
 \def \@svsec {}\else \refstepcounter {#1}\edef \@svsec%
 {\csname the#1\endcsname. \hskip .1em }\fi \@tempskipa%
 #5\relax \ifdim \@tempskipa >\z@ \begingroup #6\relax%
 \@hangfrom {\hskip #3\relax \@svsec }{\interlinepenalty%
 \@M #8.\par }\endgroup \csname #1mark\endcsname {#7}%
 \addcontentsline {toc}{#1}{\ifnum #2>\c@secnumdepth%
 \else \protect \numberline {\csname the#1\endcsname. }%
 \fi #7}\else \def \@svsechd {#6\hskip #3\@svsec #8.%
 \csname #1mark\endcsname {#7}\addcontentsline {toc}{#1}%
 {\ifnum #2>\c@secnumdepth \else \protect \numberline%
 {\csname the#1\endcsname. }\fi #7}}\fi \@xsect {#5}}

\def\section{\@startsection {section}{1}{\z@ }%
 {-3.5ex plus -1ex minus -.2ex}{2.3ex plus .2ex}{\boldface }}

\def\thebibliography#1{%
 \section *{References.\@mkboth {REFERENCES}{REFERENCES}}%
 \list {[\arabic {enumi}]}{\settowidth \labelwidth {[#1]}%
 \leftmargin \labelwidth \advance \leftmargin \labelsep %
 \usecounter {enumi}} \def \newblock %
 {\hskip .11em plus .33em minus -.07em} \sloppy \clubpenalty 4000%
 \widowpenalty 4000 \sfcode`\.=1000\relax}

\def\@maketitle{%
 \newpage \null \vskip 2em
 \begin{center}
{\Large\boldface \@title \par }
 \vskip 1.5em
 {\large \lineskip .5em
 \begin {tabular}[t]{c}\@author
 \end{tabular}\par}
 \end{center}
  \vskip .8em}

\def\abstract{%
\if@twocolumn \section *{Abstract}
 \else \small\quotation\noindent{\boldface Abstract.}\fi}

\catcode`\@=13

\font\tenbold=msbm10 scaled \magstep1
\font\sevenbold=msbm7 scaled \magstep1
\font\fivebold=msbm5 scaled \magstep1
\newfam\boldfam 
\textfont\boldfam=\tenbold
\scriptfont\boldfam=\sevenbold
\scriptscriptfont\boldfam=\fivebold

\def\ot{\otimes} 
\def\rada#1#2{#1,\ldots,#2}
\def\prada#1#2{#1 + \cdots + #2}
\def\Rada#1#2#3{#1_{#2},\dots,#1_{#3}}
\def\pa{\partial}
\def\znamenko#1{{(-1)^{#1}\cdot}}
\def\doubless#1#2{{
\def\arraystretch{.5}
\begin{array}{c}
\mbox{\scriptsize $\scriptstyle #1$}
\\
\mbox{\scriptsize $\scriptstyle #2$}
\end{array}\def\arraystretch{1}
}}
\def\qed{\hspace*{\fill}
\mbox{\hphantom{mm}$\Box$}\\}

\def\frakD{{\frak D}}
\def\PP{\mbox{\small \bf P}} \def\QQ{\mbox{\small \bf Q}}
\def\HH{\mbox{\small \bf H}} \def\TT{\mbox{\small \bf T}}
\def\ud{{\underline d}^0}
\def\malesipky{{%
{% Picture saved by xtexcad 2.4
\unitlength=.15pt
\B \hskip 2pt
\begin{picture}(80.00,40.00)(0.00,0.00)
\put(80.00,10.00){\vector(2,1){0}}
\put(80.00,30.00){\vector(2,-1){0}}
\put(36.00,20.00){\makebox(0.00,0.00){\line(0,1){20}}}
\put(44.00,20.00){\makebox(0.00,0.00){\line(0,1){20}}}
\bezier{50}(0.00,10.00)(40.00,-10.00)(80.00,10.00)
\bezier{50}(0.00,30.00)(40.00,50.00)(80.00,30.00)
\end{picture}}
\hskip 2pt \W
}}
\def\minim{{\frak M}}
\def\calD{{\cal D}}
\def\barvy{{\{\B,\W\}}}
\def\({{[ \hskip -.5mm [}}
\def\){{] \hskip -.5mm ]}}
\def\seq#1#2{\{#1\}_{{#2}}}
\def\widedef{{\hskip 4mm := \hskip 4mm}}
\def\B{{\tt B}} \def\W{{\tt W}} \def\T{{\tt T}}
\def\fb{{f_\bullet}}\def\lb{{l_\bullet}}
\def\gb{{g_\bullet}}\def\hb{{h_\bullet}}
\def\Iso{{{\cal I}so}}
\def\Span{{\rm Span}}
\def\Ker{{\rm Ker}}
\def\calG{{\cal G}}
\def\Hom{{\rm Hom}}
\def\oten#1#2{#1^{\otimes#2}}

\def\lbigbrace{{\left(\rule{0pt}{12pt}\right.}}
\def\rbigbrace{{\left.\rule{0pt}{12pt}\right)}}

\def\Riso{{\cal R}_{{\rm iso }}}
  \def\bfk{{\bf k}}
\def\F{{\cal F}}

\def\P{{\cal P}}

\def\bfV{{\bf V}}

\def\btb{{\B \to \W}}
\def\Ass{{\it Ass}}

\def\otexp#1#2{{#1^{\ot #2}}}
\def\ainft#1#2#3{{(#1,#2,#3_2,#3_3,\ldots)}}

\def\bfC{{\bf C}}
\def\bfU{{\bf U}}

\def\vlra{{\hbox{$-\hskip-1mm-\hskip-2mm\longrightarrow$}}}

\def\colorop #1(#2;#3){{#1}
   \left(\rule{0pt}{15pt}\right.
         \hskip -3mm \begin{array}{c}
	              #3\\#2
                     \end{array}
         \hskip -3mm \left. 
   \rule{0pt}{15pt} \right)
}
\def\coll#1{{\{#1(n)\}_{n\geq 1}}}
\def\id{1\!\!1}
\def\lra{{\longrightarrow}}
\def\vlra{{\hbox{$-\hskip-1mm-\hskip-2mm\longrightarrow$}}}
\def\rada#1#2{{#1,\ldots,#2}}
\def\prada#1#2{{#1+\cdots+#2}}
\def\orada#1#2{{#1\otimes\cdots\otimes#2}}
\def\otexp#1#2{{#1}^{\ot #2}}

\def\Rada#1#2#3{#1_{#2},\dots,#1_{#3}}

\def\End{\hbox{${\cal E}\hskip -.1em {\it nd}$}}
\def\bigoperad{\calF(X_\B;f,\barX; X_\W)}
\def\bfE{{\bf E}}
\def\overe{{\overline e}}
\def\underr{{\underline r}}

\def\calA{{\cal A}}
\def\calB{{\cal B}}
\def\calF{{\cal F}}
\def\frakI{{\frak I}}
\def\Asipka{{\cal A}_{\B \to \W}}
\def\Bsipka{{\cal B}_{\B \to \W}}
\def\betasipka{{\beta}_{\B \to \W}}
\def\barX{\overline X}
\def\barx{{\overline x}}

\def\colorop #1(#2;#3){{#1}
   \left(\rule{0pt}{15pt}\right.
         \hskip -3mm \begin{array}{c}
	              #2\\#3
                     \end{array}
         \hskip -3mm \left. 
   \rule{0pt}{15pt} \right)
}
\def\cases#1#2#3#4{
                  \left\{
                         \begin{array}{ll}
                           #1,\ &\mbox{#2}
                           \\
                           #3,\ &\mbox{#4}
                          \end{array}
                   \right.
}
\def\ddd{{\overline{d}_0}}

\def\calP{{\cal P}}
\def\Ainfty{{A(\infty)}}
\def\<{{< \hskip -.3mm}}
\def\susp{\uparrow\!}

\def\suspit#1{{\uparrow^{#1}\,}}

\title{HOMOTOPY DIAGRAMS OF ALGEBRAS}

\author{Martin Markl%
\thanks{The author was supported by the
grant GA AV \v CR 201/99/0675}}

\begin{document}

\maketitle

\bibliographystyle{plain}
\baselineskip 18pt plus 2pt minus 1 pt

\begin{abstract}
In~\cite{markl:ha} we proved that strongly homotopy (sh)
algebras are homotopy invariant concepts in the category of chain
complexes. 
Our arguments were based on the fact that strongly homotopy  
algebras are algebras
over minimal cofibrant operads and on the principle that
algebras over cofibrant operads are homotopy invariant. In our
approach, algebraic models for colored operads describing diagrams
of homomorphisms played an important r\^ole.

The aim of this paper is to give an explicit description of these
models. This description is necessary for
practical applications of some of the conceptual results of~\cite{markl:ha}.
Another possible application is an appropriate 
formulation of the homological perturbation lemma for chain complexes with
algebraic structures in the spirit of~\cite{markl:ip}. Our results
also provide a conceptual approach to `homotopies through
homomorphism' for strongly homotopy algebras. We also argue that
strongly homotopy 
algebras form a honest (not only weak Kan) category.

As an `application of an application' we find the homotopy
structure on the category of strongly homotopy  associative algebras 
and their strongly homotopy 
homomorphisms described in~\cite[5.4]{grandis}.

The paper can be understood as a continuation of our program to
translate~\cite{boardman-vogt:73} to algebra.
We recommend to look at Section~\ref{Misa}
first to get a broader perspective of methods, results and 
implications of this technical paper.

\vskip 2mm
\noindent 
{\boldface Subject Classification:} 55U35, 55U15, 12H05, 18G55
\hfill\break 
\noindent 
{\boldface Keywords:} colored operad, cofibrant model, homotopy
diagram
\end{abstract}

\section*{Plan of the paper. {\rm (in place of Introduction)}}

\noindent
\hangindent=5mm
\hangafter=1
{\em Section~\ref{Introduction}: Notation and terminology.\/}
We recall standard facts concerning (colored) operads and their
(minimal) models.

\noindent
\hangindent=5mm
\hangafter=1
{\em Section~\ref{prichazi_pomalu}: Minimal model for a homomorphism.\/}
We describe the minimal model $\minim_\btb$
of the colored operad $\calA_\btb$ describing a
homomorphism $f$ of $\calA$-algebras. The minimal operad
$\minim_{\btb}$ encapsulates strongly homotopy  homomorphisms of
strongly homotopy algebras.
The main result is
Theorem~\ref{Dan1}, the rest of this section is occupied by its proof.

\noindent
\hangindent=5mm
\hangafter=1
{\em Section~\ref{jak_to_vsechno_stihnu?}: Some applications\/}
of Theorem~\ref{Dan1} are given. We prove a theorem about extensions
of strongly homotopy  
algebra structures (Theorem~\ref{carman}). We also prove that
the abelization of a homotopy associative multiplication can be
extended to a balanced strongly homotopy  
associative algebra (Corollary~\ref{go}).
The remaining part of the paper
is independent on this section.

\noindent
\hangindent=5mm
\hangafter=1
{\em Section~\ref{Opulka_se_Zebrulkou}: Homotopy through homomorphisms.\/}
We construct, in Theorem~\ref{Kukacky}, another model
$\minim_\malesipky$ of $\calA_\btb$ with two independent generators
for $f$. Representations of $\minim_{\malesipky}$ describe homotopies
of homomorphisms of homotopy algebras.

\noindent
\hangindent=5mm
\hangafter=1
{\em Section~\ref{Opicak_Fuk}: Strong homotopy equivalences of
algebras.\/}
In Theorem~\ref{asi_do_te_Ameriky_neodjedu} 
we describe a model $\minim_\Iso$ 
of the colored operad $\calA_\Iso$ 
for two mutually inverse homomorphism of $\calA$-algebras.
We call colored algebras over $\minim_\Iso$ 
strong homotopy equivalences of strongly
homotopy algebras (Definition~\ref{zase_spatna_disketa}).

\noindent
\hangindent=5mm
\hangafter=1
{\em Section~\ref{Misa}: Final remarks and challenges.\/}
We discuss  possible applications and generalizations of the 
methods developed in this paper.
Namely, we discuss the `category' of strongly homotopy  
algebras, `good' homotopy
equivalences and the ideal homological perturbation lemma, and
homotopies through homomorphisms. Then we
propose a possible generalization of the main results of the paper.

\section{Notation and terminology}
\label{Introduction}

As a reference 
for standard terminology concerning operads, collections, ideals,
presentations~etc. we
recommend~\cite{may:1972,ginzburg-kapranov:DMJ94}
and~\cite{markl:zebrulka}. Algebraic models of colored operads were
studied in~\cite{markl:ha,markl:ip}.
If not said otherwise, all algebraic object will be defined over a
field ${\bf k}$ of {\em characteristic zero\/}. A map will be called a {\em
quasi-isomorphism\/} or a {\em quism\/} if it induces an isomorphism of
cohomology.

Let us recall {\em colored\/} operads describing diagrams of algebras.
Fix a (finite) {\em set of colors\/} ${\frak C}$ and
consider an operad $\P = \coll \P$ such
that each $\P(n)$ decomposes to the direct sum
\begin{equation}
\label{Kasparek}
\P(n) = \bigoplus \colorop {\P}(c;c_1,\ldots,c_n),
\end{equation}
where the summation runs over all colors $c, \Rada c1n \in {\frak C}$. 
We require the decomposition~(\ref{Kasparek}) to be, in the
obvious sense, 
$\Sigma_n$-equivariant. We also demand the following.

Let $x \in \colorop{\P}(c;\Rada c1n)$ and $x_i\in
\colorop{\P}(d_i;\rada{d^i_1}{d^i_{k_i}})$, $1\leq i \leq n$. Then we
require that the non-triviality of the 
composition $x(\Rada x1n)$ implies that
\begin{equation}
\label{Chrochtatko}
d_i = c_i, \ \mbox {for $1\leq i \leq n$}, 
\end{equation}
in which case
\[
x(\Rada x1n) \in 
\colorop{\P}(c;\rada{d^1_1}{d^1_{k_1}},\ldots,\rada{d^n_1}{d^n_{k_n}}).
\]

The intuitive meaning of~(\ref{Chrochtatko}) is that one may plug
the element $x_i$ into the
$i$-th slot of the element $x$ if and only if the color of the output
of $x_i$ is the same as the color of the $i$-th input of $x$,
otherwise the composition is defined to be zero.

An example is provided by the {\em colored endomorphism operad $\End_\bfU$\/}
on a `colored' chain complex $\bfU = \bigoplus_{c\in \bfC} \bfU_c$
where we put, in~(\ref{Kasparek}), 
\[
\colorop{\End_\bfU}(c;c_1,\ldots,c_n) := {\rm Hom}(\orada
{U_{c_1}}{U_{c_n}},U_c). 
\]
A {\em colored algebra\/} over a colored operad $\calP$ is then a map
of operads $A : \calP \to \End_\bfU$. We also sometimes say that $A$
is a {\em representation\/} of $\calP$.

\begin{example}
{\rm\
For an ordinary dg (= differential graded) 
operad $\calA = (\calA,d_\calA)$, define the $\{\B,\W\}$-colored
operad $\Asipka$ as
\begin{equation}
\label{Petrinacek1}
\Asipka := \frac{\calA_\B * \calA_\W * \calF (f)}
                  {(fa_\B = a_\W \otexp fn,\
                   \forall a\in \calA(n),\ n \geq 1)}.
\end{equation}
In this formula, $\calA_{\B}$ (resp.~$\calA_\W$) denotes the copy of $\calA$
`concentrated' in the color $\B$ (resp.~$\W$).
The symbol $f$ is a new generator, $f : \B \to \W$, and $\F(f)$
denotes the free colored operad generated by $f$. In fact, $\F(f)$
consist only of $f$, since there is no way to compose $f$ with itself.
The asterix $*$ denotes the free product of colored operads and 
the equation $fa_\B = a_\W \otexp fn$ generating the ideal in the
denominator of~(\ref{Petrinacek1}) expresses the fact that
$f$ commutes with all operations of the operad $\calA$.

There exists a differential on $\calA_\B * \calA_\W *
\calF (f)$ induced by $d_\calA$ on $\calA_\B$, $\calA_\W$ and
trivial on $\calF(f)$ which descends to a 
differential $d = d_{\Asipka}$ on $\Asipka$. 
It is clear that an algebra over $\Asipka$ consists of two
$\calA$-algebras $U$ and $V$ and a chain map $f: U \to V$ which
is also a $\calA$-homomorphism.
}
\end{example}

The following definition, introduced for ordinary operads
in~\cite{markl:zebrulka}, is crucial.

\begin{definition}
\label{minimal_model}
Let $\calA = (\calA,d_\calA)$ be an (ordinary or colored) 
dg operad. A minimal 
model of $\calA$ is a differential graded operad 
$\minim_{\calA}= (\F(E),d_{\minim})$, where ${\cal F}(E)$ is the free operad on
a collection $E$, 
together with a map $\alpha_{\calA} : \minim_{\calA} \to \calA$ of dg
operads 
such that
\begin{itemize}
\item[(i)]
$\alpha_{\calA} : \minim_{\calA} \to \calA$ is a quasi-isomorphism, and
\item[(ii)]
$d_{\minim}(E)$ consists of decomposable elements of the free operad
$\F(E)$ (the minimality).
\end{itemize}
\end{definition}

\noindent 
In~\cite[Theorem~2.1]{markl:zebrulka} 
we proved the following theorem.

\begin{theorem}
\label{e_of_mm}
Let $\calA = (\calA, d_\calA)$ be an (ordinary) operad with 
\begin{equation}
\label{Jitka}
H_*(\calA(1),d_\calA) \cong \bfk.
\end{equation}
Then there exists a minimal model $\rho
: (\calF(X),d_{\minim}) \to (\calA,d_\calA)$, unique up to isomorphism.
\end{theorem}

In fact, in~\cite[Theorem~2.1]{markl:zebrulka} we assumed that
$\calA(1)\cong {\bf k}$, but it is immediate to see that each $\calA$
satisfying~(\ref{Jitka}) is quasi-isomorphic to some
$\bar{\calA}$ with $\bar{\calA}(1) \cong {\bf k}$. The lifting
property of minimal models guarantee that the minimal model of
$\bar{\calA}$ is also a minimal model of $\calA$.

It can be immediately seen that the generators of the minimal model
from Theorem~\ref{e_of_mm} satisfy $X(1) =0$, that is, there are no
generators of arity~$1$.

As we explained in~\cite{markl:ha}, the minimal model of the
operad $\Asipka$ (unique up to an
isomorphism, by an easy modification of~\cite[Theorem~2.1]{markl:zebrulka}) 
describes {\em strongly homotopy ${\cal A}$-algebras\/} and
their {\em strongly homotopy maps\/}. This is illustrated in the following
example, taken again from~\cite{markl:ha}.

\begin{example}
\label{doufejme}
{\rm\
Let $\Ass{} := \F(\mu)/(\mu(\mu \ot \id)
- \mu(\id \ot \mu))$ be the operad for associative (non-unital) algebras.
Then the minimal
model of the operad 
\[
\Ass_\btb = \F(\mu,\nu,f)/(\mu(\mu \ot \id) = \mu(\id \ot \mu),\
\nu(\nu \ot \id) = \nu(\id \ot \nu),\ f \mu = \nu \otexp f2)
\]
is given by 
\begin{equation}
\label{Herkules}
{\Ass_\btb} \stackrel{\alpha}{\longleftarrow}
\left(
\F\lbigbrace
\mu_2,\mu_3,\mu_4,\ldots,f_1,f_2,f_3,\ldots,
\nu_2,\nu_3,\nu_4,\ldots \rbigbrace,\pa
\right),
\end{equation}
where
\begin{equation}
\label{zase_nekdo_troubi}
\def\arraystretch{1.2}
\begin{array}{l}
\mbox{
$\mu_n : \B^{\otimes n} \to \B$ is a generator of 
degree $n-2$, $n \geq 2$,}
\\
\mbox{
$f_n : \B^{\ot n} \to \W$ is a generator of 
degree $n-1$, $n\geq 1$, and}
\\
\mbox{ 
$\nu_n : \W^{\otimes n} \to \W$ is a generator of
degree $n-2$, $n \geq 2$.}
\end{array}
\def\arraystretch{1}
\end{equation}
The map $\alpha$ is defined by
 $\alpha(\mu_2) = \mu$, $\alpha(\nu_2)= \nu$, $\alpha(f_1)= f$,
while $\alpha$ is trivial on remaining generators.
The differential is given by
\begin{eqnarray}
\label{X}
\pa (\mu_n)
&:=&
\sum_{\doubless{i+j=n+1}{i,j \geq 2}} \sum_{s=0}^{n-j}
(-1)^{i+s(j+1)} \mu_{i}(\id^{\ot s} \ot \mu_j \ot
\id^{\ot i-s-1}),
\\
\label{Y}
\pa (f_n) &:=&
-\sum_{k=2}^n \sum_{\prada{r_1}{r_k} =n} 
(-1)^{\sum_{1\leq i< j \leq k}r_i(r_j+1)}
\nu_k(\orada{f_{r_1}}{f_{r_k}})+
\\
\nonumber 
&&\hskip 1cm - \hskip -3mm
\sum_{\doubless{i+j=n+1}{i\geq 1,\ j \geq 2}}
\sum_{s=0}^{n-j}
(-1)^{i + s(j+1)}
f_{i}(\id^{\ot s} \ot \mu_j \ot
\id^{\ot i-s-1}),
\\
\label{Z}
\pa (\nu_n)
&:=&
\sum_{\doubless{i+j=n+1}{i,j \geq 2}} \sum_{s=0}^{n-j}
(-1)^{i+s(j+1)} \nu_{i}(\id^{\ot s} \ot \nu_j \ot
\id^{\ot i-s-1}),
\end{eqnarray}
Let us denote the operad on the right hand side of~(\ref{Herkules}) by
$\Ainfty_\btb$. Algebras over $\Ainfty_\btb$ are clearly
$A(\infty)$-algebras and their 
$\Ainfty$-homomorphisms in the sense of~\cite{markl:JPAA92}.

The fact that the above object is indeed a minimal
model of the bi-colored operad $\Ainfty_\btb$ is a fairly
nontrivial. While it is clear that all the objects above are
well-defined, we need Proposition~\ref{c2} to show that 
$\alpha$ is a quism. See also Example~\ref{Wie}.
}
\end{example}

\section{Minimal model for a homomorphism}
\label{prichazi_pomalu}

The aim of this section is to generalize Example~\ref{doufejme} and
describe the minimal model of the $\{\B,\W\}$-colored operad 
$\Asipka$ for an arbitrary operad $\calA$ satisfying~(\ref{Jitka}) in
terms of the minimal model of $\calA$.
We need to introduce some notation first.

Let $j_{c} : \calA \to
\Asipka$ be, for $c \in \{\B,\W\}$, the map induced by the identification of
$\calA$ with $\calA_{c}\subset  \calA_\B * \calA_\W * \calF
(f)$.  It is easy to see that the map $j_{c}$ is an inclusion. 
Construction~(\ref{Petrinacek1}) is clearly functorial, that is,  
every map
$\beta: \calA  \to \calB$ induces a natural map 
$\betasipka : \Asipka \to \Bsipka$ of colored dg operads with the
property that $\betasipka(j_{c}(a)) = j_{c}(\beta (a))$, for
$a\in \calA$ and $c \in \{\B,\W\}$. 
The following proposition describes the
structure of $\Asipka$.

\begin{proposition}
\label{AaA}
Let $c, \Rada c1n \in \{\B,\W\}$ be a sequence of colors.
The component
\begin{equation}
\label{Zavodnicek}
\colorop \Asipka (c;\Rada c1n)
\end{equation}
of the operad $\Asipka$ is trivial for $c = \B$ and $(\Rada c1n) \not=
(\rada{\B}{\B})$, while it is canonically isomorphic to 
$\calA(n)$ for all other choices of colors. These canonical
isomorphisms are functorial and 
commute with the differentials.
\end{proposition}

\noindent 
{\boldface Proof.}
It is clear that~(\ref{Zavodnicek}) is trivial for $c = \B$ and 
$(\Rada c1n) \not=
(\rada{\B}{\B})$,  because there is no way to create a `map'
$(\Rada c1n) \to \B$ 
in~(\ref{Zavodnicek}) from the generators.
The only nontrivial possibility with $c = \B$
is
\[
\calA_{\B}(n) \cong \colorop \Asipka (\B;\rada {\B}{\B}),
\]
which is clearly isomorphic to $\calA(n)$. 

It follows from relations in the denominator
of~(\ref{Petrinacek1}) that, 
{}for $c = \W$, any element $a$ of~(\ref{Zavodnicek}) can be
naturally and uniquely
presented as
\[
a = \overline a (\xi(c_1) \otimes \cdots \otimes \xi(c_n)),
\]
where $\overline a \in \calA_{\W}(n) \cong \calA(n)$ and
\[
\xi(c_i) :=
\cases f{if $c_i = \B$, and}
       {\id_{\B}}{if $c_i = \W$.}
\]
Since $d(f) =0$,
the functorial correspondence $a \longleftrightarrow \overline a$
commutes with the differentials.%
\qed

Proposition~\ref{AaA} has the following important corollary.

\begin{corollary}
\label{snad_to}
Suppose that $\beta:\calA \to \calB$ is a quasi-isomorphism. Then
the induced map $\betasipka : \Asipka \to \Bsipka$ is a 
quasi-isomorphism, too.
\end{corollary}

Let $\calA = (\calA,d_\calA)$ be a dg operad satisfying~(\ref{Jitka})
and let $\rho : (\calF(X),d_{\minim}) \to (\calA,d_\calA)$
be its minimal model (see Definition~\ref{minimal_model}). 
Let $X_c$ be, for $c \in \barvy$, another copy of
the space of generators $X$ and let $i_{c} : 
X \stackrel{\cong}{\to}
X_c$ be the identification. Let $\barX := \uparrow \hskip -1mm X$ 
be the suspension of
$X$ and let $\uparrow \hskip 1mm : X \stackrel{\cong}{\to} \barX$ be the
canonical map. For $x \in X$, let $x_c := i_{c}(x)$, $c \in \barvy$, and
$\barx := \uparrow \hskip -1mm x$.

We will be working with 
the free operad $\calF(X_{\B};f,\barX; X_{\W})$
generated by the collections $X_\B$, $X_{\W}$, $\barX$ and 
a generator $f : \B \to \W$.
Generators $X_c$ will be considered `concentrated' 
in a color $c \in \barvy$ and $\barX$ as a collection of `maps' $\B
\ot \cdots \ot \B \to \W$.
We will again have canonical inclusions
$j_{c} : \calF(X) \to \calF(X_{\B};f,\barX; X_{\W})$, $c \in \barvy$. 
We will denote, for  $n\geq 1$, by 
\begin{equation}
\label{Petrinacek_a_Zahradnicek}
\frakI^{< n} \subset  \calF(X_{\B};f,\barX;
X_{\W}) 
\end{equation} 
the ideal generated by $\barX(\< n ) := \bigoplus_{k < n}\barX(k)$.
The main theorem of this section reads:

\begin{theorem}
\label{Dan1}
Let be an (ordinary) operad satisfying~(\ref{Jitka})
and let $\rho : \minim
\to (\calA,d_\calA)$, $\minim = (\calF(X),d_{\minim})$, 
be its minimal model.
Then the minimal model $\minim_\btb$ for $\Asipka$ is of the form
\begin{equation}
\label{Glue_All}
\alpha : (\calF(X_{\B};f,\barX; X_{\W}),D) \to (\Asipka,d)
\end{equation}
\label{s_Novym_godom}
such that, for each $n\geq 2$ and $x\in X(n)$,
\begin{equation}
\label{i}\label{spravny_tvar}
\def\arraystretch{1.2}
\begin{array}{cc}
      D(x_{\B}) = j_{\B}(d_{\minim}(x)),\ 
      D(x_{\W}) = j_{\W}(d_{\minim}(x)),
\\
      D(f) = 0  \mbox{ and } 
      D(\barx) = fx_\B - x_\W \otexp fn + \omega,
\end{array}
\def\arraystretch{1}
\end{equation}
for some $\omega = \omega_x \in \frakI^{< n}$ which linearly depends
on $x$.
The map $\alpha$ is given by
\begin{equation}
\label{ii}
\mbox{$\alpha(x_{\B}) = j_{\B}(\rho(x))$, 
      $\alpha(x_{\W}) = j_{\W}(\rho(x))$, 
      $\alpha(f) = f$
      and $\alpha(\barx) = 0$, \mbox { for $x \in X$.}
      }
\end{equation}
\end{theorem}

We call $fx_\B - x_\W \otexp fn$ the {\em principal part\/} and
$\omega$ {\em the tail\/} of $D(\barx)$.

\begin{remark}
{\rm\
It is well-known that if the operad $\calA$ is quadratic
Koszul~\cite{ginzburg-kapranov:DMJ94},  
its minimal model $(\calF(X),d_{\minim})$ is the cobar construction
over the quadratic dual $\calA^!$ of the operad
$\calA$~\cite[Proposition~2.6]{markl:zebrulka}. Something similar is
true also here -- it can be shown that for $\calA$ quadratic Koszul,
there exists a `closed' formula for the tail $\omega$. An example
is provided by the minimal model $A(\infty)_\btb$ of
$\Ass_{\btb}$ described in 
Example~\ref{doufejme}, compare also Example~\ref{Wie}. We are lead to:

\begin{problem}
\label{jaro}
Is there a notion of (quadratic) duality, Koszulness and of the
cobar construction for colored operads such that 
colored operad $\calA_\btb$ is Koszul if $\calA$ is and
the minimal model of $\calA_\btb$ is the cobar construction on the
dual $(\calA_\btb)^!$?
\end{problem}

Note that the differential $D$ described in Theorem~\ref{Dan1} 
is not quadratic, so one cannot expect the solution of 
Problem~\ref{jaro} to be a straightforward one. 
}
\end{remark}

The {\boldface proof of Theorem~\ref{Dan1}} 
will occupy the rest of this section. We show
first that conditions~(\ref{spravny_tvar}) and~(\ref{ii})
of Theorem~\ref{Dan1} already guarantee that
both $\alpha$ and $D$ are well defined (Lemma~\ref{c1}) and that
$\alpha$ is a quasi-isomorphism (Proposition~\ref{c2}). We then
show that the differential $D$ indeed exists. 

\begin{lemma}
\label{c1}
For any differential $D$ as in~(\ref{i}), the map $\alpha$ defined
by~(\ref{ii}) commutes with the differentials, that is
\begin{equation}
\label{Petrinacek}
d_\calA\alpha(a) = \alpha(Da),\
\forall a \in \calF(X_{\B};f,\barX; X_{\W}).
\end{equation}
\end{lemma}

\noindent
{\boldface Proof of the lemma.}
It is enough to verify~(\ref{Petrinacek}) on generators. It obviously
holds  for $a\in X_{\B}$, $a\in X_{\W}$ or $a=f$. 
For $a= \barx \in \barX$, equation~(\ref{Petrinacek}) reduces to
$0 = \alpha (D(\barx))$. We have
\[
\alpha (D(\barx)) = \alpha(fx_\B - x_\W \otexp fn + \omega)
= f j_{\B}(\rho(x)) - j_{\W}(\rho(x)) \otexp fn + \alpha(\omega),
\]
because $\alpha$ is a homomorphism. It follows from the definition of
$\Asipka$ that
$f j_{\B}(\rho(x)) - j_{\W}(\rho(x)) \otexp fn = 0$. Moreover, 
$\alpha(\omega) = 0$, since
$\omega \in \frakI^{< n}$ and 
$\alpha$ is trivial on the generators $\barX(<n )$ of $\frakI^{< n}$.
Thus $\alpha(D(\barx)) = 0$ which finishes the proof.%
\qed

\begin{proposition}
\label{c2}
Suppose that the differential $D$ is as in~(\ref{i}).
Then the map $\alpha$ defined by~(\ref{ii}) is a 
quasi-isomorphism.
\end{proposition}

\begin{example}
\label{Wie}
{\rm\
The minimal model of the operad $\Ass$ for associative algebras is
$(\calF(X),d_{\minim})$ with $X$ generated by $\mu_2,\mu_3,\ldots$ as
in~(\ref{zase_nekdo_troubi}) and the
differential $d_{\minim}$ given by~(\ref{X}). Equation~(\ref{Y}) is clearly
of the form
\[
\pa f_n =  f_1 \mu_n - \nu_n f_1^{\otimes n} + \omega
\]
where $\omega$ belongs to the ideal generated by $f_2,f_3,\ldots,f_{n-1}$.
Therefore the model for $\Ass_{\btb}$
described in Example~\ref{doufejme} is of the type predicted by
Theorem~\ref{Dan1}, with $X_\B = \Span(\mu_2,\mu_3,\ldots)$, $X_\W
= \Span(\nu_2,\nu_3,\ldots)$ and $\barX = \Span(f_2,f_3,\ldots)$. In
particular, the map $\alpha$ constructed in this example is a 
quasi-isomorphism.
}
\end{example}

To prove Proposition~\ref{c2}, observe that the map $\alpha$
of~(\ref{Glue_All}) decomposes as
\begin{equation}
\label{decomposition}
(\calF(X_{\B};f,\barX; X_{\W}),D) \stackrel{\gamma}{\lra}
  \left(
      \frac{\calF(X_{\B};f;X_{\W})}
                        {(fx_{\B} - x_{\W} \otexp fn)}, d 
  \right) =\minim_\btb
\stackrel{\rho_\btb}{\vlra}
\calA_\btb.
\end{equation}
Since the map $\rho_\btb$ is a quism by Corollary~\ref{snad_to}, it is
enough to prove that the map $\gamma$ is also a quism:

\begin{lemma}
\label{32}
For each minimal operad $(\calF(X),d_{\minim})$ is the map
\[
\gamma : (\bigoperad,D) \longrightarrow (\calF(X),d_{\minim})_{\B\to\W}
\]
given by $\gamma(x_{\B}) = j_{\B}(x)$, $\gamma(x_{\W}) = j_{\W}(x)$,
$\gamma(f) =f$  and
$\gamma(\barx) = 0$, a quasi-isomorphism.
\end{lemma}

\noindent
{\boldface Proof the lemma} is based on a repeated spectral sequence argument.
Let us define an ascending filtration of $\calF(X_\B;f,\barX; X_\W)$,
\begin{equation}
\label{napred_sluzba_potom_druzba}
{\frak F}_1: \hskip .5cm
\cdots \subseteq F_{-2}  \subseteq F_{-1} \subseteq F_{0} = F_1 = F_2
= \cdots = \bigoperad,
\end{equation}
by postulating $F_p$ to be the subspace spanned by 
expressions having at least $-p$ occurrences of generators from $X_{\B},
X_{\W}$ or $\barX$. Since
clearly $D(F_p) \subset F_p$, we can consider the related
spectral sequence $\bfE = (E^r_{pq},d^r)$.
This spectral sequence converges because, for any fixed arity
$n$, $F_p(n)= 0$ for $p$ sufficiently small.
The initial term of $\bfE$ is easy to describe, 
\[
(E^0,d^0) \cong (\bigoperad,d^0),
\]
with $d^0$ denoting the `linear part' of $d$, namely
\begin{equation}
\label{transfer}
d^0(x_{\B}) = d^0(x_{\W}) = d^0(f) =0 \mbox { and }
d^0(\barx) = fx_{\B} - x_{\W} \otexp fn, \mbox { for } x \in X(n).
\end{equation}
We claim that
\begin{equation}
\label{44}
H_*(E^0,d^0) \cong \frac{\calF(X_{\B};f;X_{\W})}
                        {(fx_{\B} - x_{\W} \otexp fn)}.
\end{equation}                    
Let us denote, just for the purpose of this proof, $\calF :=
\calF(X_{\B};f,\barX;X_{\W})$. The colored operad $\calF$ has also another
`upper' grading, $\calF = \bigoplus_{s\geq 0}{\calF}^s$, 
induced by the number of generators from
$\barX$. The differential $d^0$ lowers this upper degree by $1$.
It is immediate to see that
\[
H_*^0(\calF^*_*,d^0) 
\cong \frac{\calF(X_{\B};f;X_{\W})}
                        {(fx_{\B} - x_{\W} \otexp fn)},
\]
equation~(\ref{44}) will thus follow from $H_*^{> 0}(\calF^*_*,d^0) = 0$. To
prove this, 
observe that there is yet another ascending filtration of $\calF$ induced by
the number of occurrences of the generator $f$, i.e.~$\calF = \bigcup
G_p$, where $G_p$ is spanned by expressions having at least $-p$
occurrences of the generator $f$.
The first term of the
related spectral sequence is $(\calF^*_*,\ddd)$, where
\[
\ddd(x_{\B}) = \ddd(x_{\W}) = \ddd(f) =0 \mbox { and } \ddd(\barx) =
fx_{\B},
\mbox { for $x \in X$}.
\]
We shall prove that $(\calF^*_*,\ddd)$ is acyclic in positive `upper'
gradings. This means that 
\begin{equation}
\label{22}
\left( \colorop {\calF^{> 0}_*} (c; \Rada c1n), \ddd \right)
\end{equation}
is acyclic for each choice of colors $c, \Rada c1n \in
\{\B,\W\}$. 
If $c = \B$
then, as in the proof of Proposition~\ref{AaA}, 
we conclude that~(\ref{22}) cannot
contain a generator from $\barX$, thus it is not only acyclic, but
even trivial in positive
upper dimensions.

If $c = \W$, each element $a$ of~(\ref{22}) can be expanded as 
\begin{equation}
\label{deprese_z_Rataj}
a = \sum_{m \geq 2}(f e^{m,i}_\B r^{\B}_{m,i}
+ \overe^{m,i}\underr_{m,i} + e^{m,i}_\W r^{\W}_{m,i}
)
\end{equation}
where $\{e^{m,i}\}$ is a basis of $X(m)$ and
$r^{\B}_{m,i},r^{\W}_{m,i},
\underr_{m,i} \in \otexp \calF m$, $m \geq 2$ (the summation
convention assumed).
We prove, by induction on the arity $n$, that  
\begin{equation}
\label{zelena}
\mbox {
if $a$ is of a positive upper grading and $\ddd(a) = 0$, then $a$
is a $\ddd$-boundary.}
\end{equation}
For $n=2$, equation~(\ref{deprese_z_Rataj}) reduces to
$a = \sum_i  \underline{\alpha}_i \overe^{2,i}$
for some $\underline{\alpha}_i \in \bfk$, thus 
$\ddd(a) = \sum_i \underline{\alpha}_i f e_\B^{2,i} = f (\sum_i 
\underline{\alpha}_i e_\B^{2,i})$.
The condition  $\ddd(a) = 0$ then
implies that $\underline {\alpha}_i = 0$
for each $i$, thus $a =0$.

Suppose we have proved~(\ref{zelena}) 
for all arities $< n$ and suppose that $a$ has
arity $n$.
If $\ddd(a) = 0$, then
\[
0 =
\sum_{m \geq 2}(
\znamenko{e^{m,i}} f e^{m,i}_\B  \ddd(r^{\B}_{m,i})
+ f e^{m,i}_\B\underr_{m,i} 
-\znamenko{e^{m,i}} \overe^{m,i}\ddd(\underr_{m,i})
+ \znamenko{e^{m,i}} e^{m,i}_\W \ddd(r^{\W}_{m,i})
),
\]
which implies that, for all $m \geq 2$, 
\[
\ddd(\underr_{m,i}) = 0,\
\ddd(r^\B_{m,i}) = - \znamenko{e^{m,i}}\underr_{m,i}
\mbox { and }
\ddd(r^\W_{m,i}) = 0.
\]
Since each $r^\W_{m,i}$ is a product of elements of arity $< n$, by
induction assumption there exists $b^\W_{m,i}$ such that
$r^\W_{m,i} = \ddd(b^\W_{m,i})$.
It is easy to verify that then
\[
a = \ddd \left\{
      \sum_{m\geq 2}(\overe^{m,i}r^{\B}_{m,i} + 
      (-1)^{e^{m,i}} e_\W^{m,i} b^\W_{m,i})
         \right\},
\]
which proves~(\ref{zelena}) for arity $n$.
Thus $(\calF^*_*,\ddd)$ 
is acyclic in positive upper dimensions, and so
is, by a spectral sequence argument, also $(\calF^*_*,d^0)$. This
proves~(\ref{44}).

Define an ascending filtration ${\frak F}_2$ of the operad
\[
(\calF(X),d_{\minim})_\btb =
\left(
      \frac{\calF(X_{\B};f;X_{\W})}
                        {(fx_{\B} - x_{\W} \otexp fn)},
      d 
\right),
\]
by the number of generators from $X_\B$ and $X_\W$. 
The differential induced on the first term of the related spectral
sequence is trivial. The map $\gamma$ is clearly a morphism of ${\frak
F}_1$-${\frak F}_2$ 
filtered operads. 
It follows
from~(\ref{44}) that $\gamma$ induces a quism of the initial terms of the
spectral sequences, thus
$\gamma$ is a quism as well. This finishes the proof of the lemma.%
\qed

To finish the proof of Theorem~\ref{Dan1}, we will need also a
`restricted' version of Lemma~\ref{32}:

\begin{lemma}
\label{restricted}
Let $K \geq 2$ and let $d_{\minim}^{<K}$ be the restriction of $d_{\minim}$ to
$X(\<K)$. 
Let $D^{<K}$ be the similar obvious restriction of the
differential $D$. Then the map
\begin{equation}
\label{56}
\gamma_{<K} 
: (\calF(X_{\B}(\< K);f,\barX(\< K);X_{\W}(\< K)),D^{<K}) 
\longrightarrow (\calF(X(\<K)),d_{\minim}^{<K})_{\B\to\W}
\end{equation}
given by 
$\gamma_{<K}(x_{\B}) = j_{\B}(x)$, $\gamma_{<K}(x_{\W}) = j_{\W}(x)$,
$\gamma_{<K}(f)=0$ and
$\gamma_{<K}(\barx) = 0$ is a quasi-isomorphism.
\end{lemma}

\noindent 
{\boldface Proof of the lemma.}
Since all objects in the lemma are well-defined and the
restricted differential $D^{<K}$ has the form~(\ref{spravny_tvar}), the map
$\gamma_{<K}$ is a quism by Lemma~\ref{32} applied to the minimal
operad $(\calF(X(\<K)),d_{\minim}^{<K})$.%
\qed

\noindent
{\boldface Proof of Theorem~\ref{Dan1}.}
In the light of Lemma~\ref{c1} and Proposition~\ref{c2}, 
it suffices to construct a
differential $D$ satisfying~(\ref{i}). We proceed by
induction on the arity $n$. 
For $\barx \in \barX(2)$, let $D(\barx) := f x_{\B} - x_{\W} \otexp f2$.
Suppose we have already defined $D$ on $\barX(< \hskip -1.5mm n)$. 
If we put, for $\barx\in
\barX(n)$ and for some $\omega \in \frakI^{< n}$,
\[
D(\barx): = fx_\B - x_\W \otexp f{n} + \omega,
\]
then $(D \circ D)(\barx) =0$ implies that $\omega$ must satisfy
\begin{equation}
\label{lhs}
D(x_{\W}) \otexp f{n} - f D(x_{\B}) = D(\omega).
\end{equation}
Let us denote the left hand side of the above equation by $\varphi$,
\[
\varphi:= D(x_{\W}) \otexp f{n} - f D(x_{\B})
= j_{\W}(d_{\minim}(x)) \otexp f{n} - f j_{\B}(d_{\minim}(x)).
\] 
We need to solve the equation
\begin{equation}
\label{hlinena_dymka}
\varphi = D(\omega)
\end{equation}
with some $\omega \in {\frak I}^{<n}$.
We are going to apply Lemma~\ref{restricted}.
To simplify the notation, denote 
\[
\calG_{< n} :=  (\calF(X_{\B}(\< n);f,\barX(\<
n);X_{\W}(\< n)),D^{<n}).
\]
Observe that ${\frak I}^{<n} \subset \calG_{< n}$.
It is also evident that $\varphi \in \calG_{< n}$ and that $D^{<
n}(\varphi) = 0$.
Since $\phi \in {\rm Ker}(\gamma_{< n})$, by
Lemma~\ref{restricted} with $K=n$
there exits $\alpha_1 \in \calG_{< n}$ such
that $\varphi = D^{< n} {\alpha_1}$. Clearly $\alpha_1$ decomposes as 
$\alpha_1 = u_1 + \omega_1$,
where $\omega_1 \in {\frak I}^{< n}$ and $u_1$ does not contain
generators from $\barX$. 

Let us apply the map $\gamma_{<n}$ to the equation 
$\varphi = D^{<n} u_1 + D^{<n} \omega_1$. 
Since $\gamma_{<n}(\varphi) = \gamma_{<n}(D^{<n}\omega_1) = 0$,
$D^{<n}u_1 
\in \Ker(\gamma_{<n})$ as well. Because $u_1$ does not contain elements of
$\barX$, $D^{<n} u_1 \in {\cal G}_{< n-1}$ thus, in fact, 
$D^{<n} u_1 \in {\rm Ker}(\gamma_{< n-1})$. 

Lemma~\ref{restricted} with $K=n-1$ gives some $\omega_2 
\in {\frak I}^{< n-1}$ and some $u_2 \in \calG_{< n-1}$ which
does not contain  generators from $\barX$
such that 
$D^{< n} u_1 = D^{< n-1} u_2 + D^{< n-1} \omega_2$. 
Repeating this process $n-2$ times, we end up with 
$D^{<4}u_{n-3} = D^{<3}u_{n-2} + D^{<3}\omega_{n-2}$, where
$\omega_{n-2} \in {\frak I}^{<3}$ and $u_{n-2} \in \calG_{<3}$  
does not contain generators from $\barX$. Thus $D^{<
3}u_{n-2} =0$ and $\omega
:= \omega_1 + \omega_2 + \ldots + \omega_{n-2}$ solves~(\ref{hlinena_dymka}).%
\qed

\section{Some applications}
\label{jak_to_vsechno_stihnu?}

In this section we present a couple of applications of the model
studied in Section~\ref{prichazi_pomalu}. Though the
results are formulated for very specific examples of homotopy
algebras, as we will see below, 
only the principal parts of the differential matter. 
Therefore similar statements can be 
formulated for any type of strongly homotopy algebras.

Recall that, for each $K \geq 1$, an {\em $A(K)$-algebra\/} is an object ${\bf
W} = (W,\pa, n_2,n_3,\ldots,n_K)$ such that the multilinear maps
$n_i : \otexp Wi \to W$, $2 \leq i \leq K$, 
satisfy all axioms of $A(\infty)$-algebras in
arities $\leq K$. Similarly, a {\em morphism\/} ${\bf F}$ of two
$A(K)$-algebras, 
\[
{\bf F} :{\bf V} = (V,\partial, m_2,m_3,\ldots,m_K) \longrightarrow 
{\bf W} = (W,\pa, n_2, n_3, \ldots,n_K)
\] 
is a sequence of multilinear maps $F_i :
\otexp Vi \to W$, $1 \leq i \leq K$, 
that satisfies all axioms of a homomorphism of
$A(\infty)$-algebras in arities $\leq K$. This notion
slightly differs from the one of~\cite[p.~147]{markl:JPAA92}.

Each $A(\infty)$-algebra ${\bf V}$ determines, by forgetting all structure
operations of arities $> K$, an $A(K)$-algebra ${\bf V}_K$.

\begin{theorem}
\label{carman}
Suppose we are given, for some fixed $K \geq 1$, the following data:
\begin{itemize}
\item[(i)] An $A(\infty)$-algebra ${\bf V} =
(V,\partial, m_2,m_3,\ldots)$,
\item[(ii)]
an $A(K)$-algebra ${\bf W} = (W,\partial, n_2,n_3,\ldots,n_K)$ and
\item[(iii)]
a morphism $\{F_n: V ^{\otimes i} \to W\}_{n\leq
K} : {\bf V}_K \to {\bf W}$ of $A(K)$-algebras. 
\end{itemize} 
Suppose that the chain map $F_1:(V,\partial) \to
(W,\partial)$ is a quasi-isomorphism. Then the above data can be
extended to an $A(\infty)$-structure ${\bf W}$ on $(W,\partial)$ and to
an $A(\infty)$-homomorphism ${\bf F}:{\bf V} \to {\bf W}$. 
A similar statement holds also for balanced sh{} associative ($C(\infty)$)
and sh{} Lie ($L(\infty)$) algebras (see~\cite[p.~148]{markl:JPAA92} 
and~\cite{lada-stasheff:IJTP93} for the
definitions of these objects).
\end{theorem}
  
\noindent 
{\boldface Proof.}
Let us assume the notation of Example~\ref{doufejme}. The data of the
theorem can be encoded by a map
\[
\phi: 
\calF(\mu_2, \mu_3, \ldots; f_1,f_2,\ldots,f_K;
\nu_2,\nu_3,\ldots,\nu_K) \to \End_{V,W}
\]
with $\phi(\mu_i) := m_i$, $2 \leq i$, $\phi(\nu_j): = n_j$, $2 \leq j
\leq K$ and $\phi(f_i): = F_i$, $1 \leq i \leq K$.
We need some
$n_{K+1} \in \Hom(\oten W{(K+1)},W)$ and 
$F_{K+1} \in \Hom(\oten V{(K+1)},W)$ such that
\begin{eqnarray}
\label{B}
\pa n_{K+1} & =& \phi(\pa \nu_{K+1}) \mbox { and}
\\
\label{A}
\pa F_{K+1} &=& 
\phi(\nu_{K+1}\oten {f_1}{(K+1)} - f_1 \mu_{K+1} + \omega),
\end{eqnarray}
$\phi(\nu_{K+1}) := n_{K+1}$ and $\phi(f_{K+1}) := F_{K+1}$  
will then extend our data one step up.
We start by~(\ref{B}). Applying $\phi$ to
\[
0 = \pa^2 f_{K+1} = 
\pa \nu_{K+1} f_1^{\otimes (K+1)} -  f_1 \pa \mu_{K+1} + \pa \omega
\]
we obtain
\[
\phi(\pa \nu_{K+1})F_1^{\otimes (K+1)} = F_1 \phi (\pa \mu_{K+1}) -
\phi(\pa\omega) = \pa(F_1\phi(\mu_{K+1} - \omega)),
\]
so $\phi(\pa \mu_{K+1})F_1^{\otimes (K+1)}$ is homologous 
to zero in $\Hom(\otexp V{(K+1)},W)$ .
Since $F_1$ is a quism,
$\phi(\pa \nu_{K+1})$ is homologous to zero in $\Hom(\otexp
W{(K+1)},W)$ and the existence of $n_{K+1}$ satisfying~(\ref{B})
easily follows.

We solve~(\ref{A}) invoking a very useful trick of an `additive
renormalization.' 
Observe that
if $\theta_{K+1} \in \Hom(\otexp W{(K+1)},V)$ is closed, 
then changing $n_{K+1}$ to $n_{K+1} +
\theta_{K+1}$ does not violate~(\ref{B}). For such a renormalized
$n_{K+1}$, the right hand side of~(\ref{A}) reads
\[
\chi := 
n_{K+1}F_1^{\otimes (K+1)} - F_1 m_{K+1}
+ \phi(\omega) + \theta_{K+1} F_1^{\otimes {(K+1)}}.
\]
The first three terms of $\chi$ are closed. Since $F_1$ a quism, there 
exists $\theta_{K+1}$ such that $\chi$ homologous to zero, and the existence
of $F_{K+1}$ solving~(\ref{A}) follows. The induction may go on.%
\qed

For $K=1$ the data of Theorem~\ref{carman}
consists of an $A(\infty)$-algebra $\bfV =
\ainft V\partial m$ and a chain equivalence $F_1 :(V,\pa) \to
(W,\pa)$. For $K=2$ we have also a bilinear product $n_2 : W\ot W
\to W$ and a homotopy $F_2 : V\ot V \to W$ between $F_1 m_2$ and
$n_2(F_1,F_1)$. 

A statement similar to Theorem~\ref{carman} can be clearly
formulated for any type of strongly homotopy algebras, since all we
needed was that the minimal model was of the form of Theorem~\ref{Dan1}.

\begin{corollary}
\label{go}
Suppose that $(U,\mu,\pa)$ is a differential graded associative
algebra and let $\nu : U\ot U \to U$ be a bilinear multiplication
chain homotopic to $\mu$. Then $\nu$ can be extended to a strongly
homotopically associative 
structure on $(U,\pa)$.   
\end{corollary}

\noindent
{\boldface Proof.} The corollary immediately follows from
Theorem~\ref{carman} with $\bfV = (V,\pa,\mu,0,0\ldots)$, $(W,\pa): =
(U,\pa)$, $F_1 := {\rm id}$ and $F_2$ the homotopy between $\mu$ and
$\nu$.%
\qed

As observed by Jim Stasheff, the above corollary must be
`spiritually' true if one believes that $A(\infty)$-structures are
homotopically invariant versions of associative algebras.
And indeed, it is very easy to prove that the multiplication $\nu$,
homotopic to an associative one, is associative up to a homotopy; it
is even possible to give an explicit formula for the homotopy. But it
is not clear whether this homotopy  extends to a coherent
hierarchy of higher homotopies, and this the hard part of
Corollary~\ref{go}.
In the following corollary, the characteristic $0$ assumption is very crucial. 

\begin{corollary}
Suppose that $(U,\pa,\mu)$ is an associative, homotopy commutative
algebra. Let $\overline \mu(u,v) := \frac 12 (\mu(u,v) + \mu(v,u))$ be
the symmetrization. Then $\overline \mu$ can be extended to a balanced
$A(\infty)$-structure on $(U,\pa)$.
\end{corollary}  
 
\noindent
{\boldface Proof.} Let $H := H_*(U,\pa)$ and let $*$ be the multiplication
induced by $\mu$ on $H$. It is immediate to see that it is
commutative associative and that it coincides with that induced by 
$\overline \mu$. Now choose a monomorphism $(H,0)
\stackrel{\iota}{\to} (U,\pa)$ inducing the identity map on
homology. It can always be done, since we are in characteristic
zero. Then $\iota (*) : H\ot H \to U$ and $\overline \mu(\iota,\iota):
H\ot H \to U$ induce the same map in cohomology, thus (again the
characteristics zero argument) there exists a homotopy $h$ between
$\iota( *)$ and $\overline \mu(\iota,\iota)$. 

The corollary now follows from Theorem~\ref{carman} with $\bfV :=
(H,0,*,0,\ldots)$, $(W,\pa,\nu_2) := (U,\pa,\overline \mu)$, $F_1 =
\iota$ and $F_2 := h$.\qed

\section{Homotopy through homomorphisms}
\label{Opulka_se_Zebrulkou}

In Section~\ref{prichazi_pomalu} we described the minimal model for the operad
$\calA_\btb$ and observed that it describes strongly homotopy
homomorphisms of strongly homotopy algebras. The aim of this section
is to understand homotopies between these homomorphisms. By our
philosophy, we need to resolve the operad
\[
\calA_{\malesipky}:= 
\frac{\calA_\B * \calA_\W * \calF (p,q)}
                  {(pa_\B = a_\W \otexp pn,\ qa_\B = a_\W \otexp qn,\ p=q)}.
\]
where $p,q: \B \to \W$ are generators of degree $0$ and 
the first two relations in the denominator
are satisfied for all $a\in \calA(n)$ and $n
\geq 1$. This operad describes two identical homomorphisms of
$\calA$-algebras and its resolution should replace the strict
equality $p=q$ by a homotopy. 

Since the operad $\calA_{\malesipky}$  is
clearly isomorphic to $\calA_\btb$, its minimal model will not give
what we want and we shall consider a resolution which has two different
generators for the same map $p=q$ instead. This resolution will not be
minimal, but it will still be cofibrant in a suitable
sense~\cite{markl:ha}. 
A toy model for this resolution is the $\barvy$-colored operad
\begin{equation}
\label{ale_prijde}
\calD = (\calF(p,q,h),d),\
d(p) = d(q):= 0 \mbox { and }
d(h): = p-q
\end{equation}
($\calD$ from dull), 
where $p,q : \B \to \W$ are generators of degree $0$ and $h : \B \to
\W$ has degree $1$. Operad~(\ref{ale_prijde}) clearly resolves
the free operad $(\calF(f),d =0)$ on one degree $0$ generator $f: \B \to
\W$. An algebra over $\calD$ consists of two chain maps and a chain
homotopy between them.
For each $m \geq 1$ the space  $\calD(\otexp{\B}m, \otexp{\W}m)$ 
contains a special degree $1$ morphism $\(h\)_{{\rm ns}}$ defined by
\[
\(h\)_{{\rm ns}} := h \otimes\otexp q{m-1} + 
                      p \otimes h \otimes \otexp q{m-2} +
                      \cdots + \otexp p{m-2} \otimes h \otimes q +
                      \otexp p{m-1}  \otimes h 
\]
and, since ${\rm char}(\bfk)=0$, also its symmetrization
\[
\(h\) := \frac 1{m!} \sum_{\sigma \in \Sigma_m} \(h\) \sigma
\]
($\Sigma_m$ is the symmetric group on $m$ symbols).
The above morphisms satisfy
the following differential equation in $\calD$
\[
d(\(h\)_{{\rm ns}}) = d(\(h\)) = \otexp p m - \otexp q n
\]
and the morphism $\(h\)$ is $\Sigma_n$-equivariant. The
maps $\(h\)$ and $\(h\)_{\rm ns}$ 
are very simple examples of a polarization which we
introduce in Definition~\ref{snezenky}.

Let $X$ be a collection. We will work with the free
operad $\calF(X_\B; p,q,h; X^p,X^q,X^h; X_\W)$ generated by
\begin{itemize}
\item[(i)]
two copies  $X_\B$ and $X_\W$ of $X$ interpreted as in
Theorem~\ref{Dan1}, 
\item[(ii)] 
the generators $p,q,h$ as in~(\ref{ale_prijde}),
\item[(iii)]
two copies of $X^p$ and $X^g$ of the suspension $\susp X$ interpreted
as collections of multilinear maps from $\B$ to $\W$ and
\item[(ii)]
a copy $X^h$ of the double suspension $\suspit2X$, again interpreted
as a collection of multilinear maps from $\B$ to $\W$. 
\end{itemize}
The notation used below is paralell to that
of Theorem~\ref{Dan1} and we believe its meaning is clear.

\begin{theorem}
\label{Kukacky}
Let $\calA$ be an (ordinary) operad satisfying~(\ref{Jitka}) and
let $\rho : \minim 
\to (\calA,d_\calA)$, $\minim =(\calF(X),d_{\minim})$, 
be the minimal model of $\calA$. Then there
exists a cofibrant model $\minim_{\malesipky}$
of the colored operad $\calA_\btb$ of the form
\[
\psi:
(\calF(X_\B; p,q,h; X^p,X^q,X^h; X_\W),D) \to (\calA_\btb,d)
\]
where the tails $\psi(x_\B) = j_\B(\rho(x))$ and $\psi(x_\W) = j_\W(\rho(x))$
for $x \in X$,
$\psi(p) = \psi(q) = f$ and $\psi$ is trivial on remaining generators.
The differential is given by
$D(p) = D(q) = 0$,
$D(h) = p-q$ and, for $x \in X(n)$, $n \geq 2$,
$D(x_\B) = j_\B(d_{\minim}(x))$, $D(x_\W) = j_\W(d_{\minim}(x))$. Moreover,
\[
\def\arraystretch{1.25}
\begin{array}{c}
D(x^p) =  px_\B - x_\W \otexp pn + \omega_{p},\ 
D(x^q) = qx_\B - x_\W \otexp qn + \omega_q, \mbox { and}
\\
D(x^h) = x^p -x^q - h x_\B + (-1)^x x_\W \(h\) + \omega_h, 
\end{array}
\def\arraystretch{1}
\]
where the tails 
$\omega_p$, $\omega_q$ and $\omega_h$ depend linearly on $x$ and 
\begin{itemize}
\item[(i)]
$\omega_p$ belongs to the ideal of the suboperad
$\calF(x_\B;p,X^p;x_\W)$ generated by $X^p(\<n)$,
\item[(ii)]
$\omega_q$ belongs to the ideal of the suboperad
$\calF(x_\B;q,X^q;x_\W)$ generated by $X^q(\<n)$, and
\item[(iii)]
$\omega_h$ belongs to the ideal generated by $X^p(\<n), X^q(\<n)$ and
$X^h(\<n)$.
\end{itemize}
\end{theorem}

\noindent
{\boldface Proof.}
The tails $\omega_p$ and $\omega_q$ can be constructed as follows.
Let $(\calF(X_{\B};f,\barX; X_{\W}),D)$ be as in Theorem~\ref{Dan1} and
consider two homomorphisms
\begin{equation}
\label{jsem_zvedav}
{% Picture saved by xtexcad 2.4
\unitlength=.7pt
\begin{picture}(40.00,25.00)(0.00,12.00)
\put(0.00,0.00){\makebox(0.00,0.00)[r]%
                     {$\theta_q: \calF(X_{\B};f,\barX; X_{\W})$}}
\put(0.00,30.00){\makebox(0.00,0.00)[r]%
                     {$\theta_p: \calF(X_{\B};f,\barX; X_{\W}) $}}
\put(40.00,10.00){\makebox(0.00,0.00)[lb]%
                     {$\calF(X_\B; p,q,h; X^p,X^q,X^h; X_\W)$}}
\put(10.00,0.00){\vector(2,1){20.00}}
\put(10.00,30.00){\vector(2,-1){20.00}}
\end{picture}}
\end{equation}
defined by 
\begin{eqnarray*}
&\theta_p(f) = p,\ \theta_p(x_\B) = \theta_q(x_\B) := x_\B,\ 
\theta_q(f) = q,&
\\
& \theta_p(x_\W) = \theta_q(x_\W) := x_\W,\ 
\theta_p(\barx) := x^p \mbox { and } \theta_q(\barx) := x^q&
\end{eqnarray*}
for $x \in X$. We put then
$\omega_p := \theta_p(\omega) \mbox { and } \omega_q :=
\theta_q(\omega)$, 
where $\omega$ is as in~(\ref{spravny_tvar}). The construction
of the differential $D$ will be completed by proving the existence of some
$\omega_h$ belonging to the ideal defined in (iii) and satisfying
\begin{equation}
\label{ten_caj_je_dobry}
D(\omega_h) = D(\omega_q - \omega_p).
\end{equation}
Assume, as in the proof of Theorem~\ref{Dan1}, 
that the differential
$D$ exists and prove that the map $\psi: \minim_{\malesipky} \to
\calA_{\btb}$ 
is a quism. Since $\alpha = \psi \theta_p = \psi \theta_q$, where
$\alpha$ is the quism of~(\ref{s_Novym_godom}), it suffices to prove
that $\theta_p$ (or $\theta_q$) is a quism.

Let ${\frak F}_3$ be the filtration of $\minim_{\malesipky}$ given by 
the number of generators from $X^p$, $X^q$ and $X^h$.
The differential $d^0$ of the first term of the 
corresponding spectral sequence 
is given by 
\begin{eqnarray*}
& d^0(h) = p-q,\ d^0(x_\B) = d^0(x_\W) = 0,\ d^0(p) = d^0(q) = 0,&
\\
&d^0(x^p) = px_\B - x_\W \otexp pn,\
d^0(x^q) = qx_\B - x_\W \otexp qn \mbox { and }&
\\
&d^0(x^h) = x^p -x^q - h x_\B  +(-1)^x x_\W \(h\),&
\end{eqnarray*}
for $x \in X$.
The substitution 
\begin{eqnarray*}
&p \mapsto p,\  h \mapsto h,\ q \mapsto -q+p,\  
x_\B \mapsto x_B,\  x_\W \mapsto x_\W,&
\\
&x^p \mapsto x^p,\  x^h \mapsto x^h \mbox { and } 
x^q \mapsto x^p - x^q - h x_\B +(-1)^x x_\W \(h\)&
\end{eqnarray*}
converts $d^0$ to the differential $\ud$ with
\begin{eqnarray*}
&\ud(p) = \ud(q) = 0,\ \ud(h) = q,\ \ud(x_\B) = \ud(x_\W) = 0,&
\\
&\ud(x^p) = px_\B - x_\W \otexp pn,\
\ud(x^q) = 0 \mbox { and } \ud(x^h) = x^q.&
\end{eqnarray*}
The dg operad $(\calF(X_\B; p,q,h; X^p,X^q,X^h; X_\W),
\ud)$ is clearly isomorphic to the free product
\begin{equation}
\label{Godom}
(\calF(X_{\B};f,\barX; X_{\W}),d^0) * (\calF(q,h,;X^q,X^h),\ dh=q,\ dx^h=x^q),
\end{equation}
where the differential $d^0$ of the first factor is as
in~(\ref{transfer}) and
the second factor is acyclic. 

Recall that in~(\ref{napred_sluzba_potom_druzba})
we defined filtration ${\frak F}_1$ of $(\calF(X_{\B};f,\barX;
X_{\W}),D)$. 
The map $\theta_p$ is clearly an ${\frak F}_1$-${\frak F}_3$ 
homomorphism of filtered operads. On the initial stages 
of the related spectral sequences 
it induces precisely the inclusion of the first factor
of~(\ref{Godom}), which is a quism. 
Thus $\theta_p$ is a quism, by a standard
spectral sequence argument.

The existence the tail $\omega_h$ satisfying~(\ref{ten_caj_je_dobry}) 
follows, as in the proof of Theorem~\ref{Dan1}, 
from a general nonsense. We leave the details to the reader.%
\qed

In the proof of Theorem~\ref{Kukacky} we in fact proved:

\begin{corollary}
The maps $\theta_p$ and $\theta_q$ of~(\ref{jsem_zvedav})
fit to the following commutative
diagram of quasi-isomorphisms of colored operads:
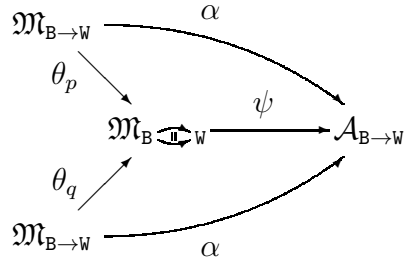
\begin{center}
{% Picture saved by xtexcad 2.4
\unitlength=1.000000pt
\begin{picture}(130.00,85.00)(0.00,0.00)
\put(110.00,30.00){\vector(1,1){0.00}}
\put(110.00,50.00){\vector(1,-1){0.00}}
\bezier{200}(20.00,0.00)(70.00,0.00)(110.00,30.00)
\bezier{200}(20.00,80.00)(70.00,80.00)(110.00,50.00)
\put(80.00,50.00){\makebox(0.00,0.00){$\psi$}}
\put(10.00,20.00){\makebox(0.00,0.00)[r]{$\theta_q$}}
\put(60.00,85.00){\makebox(0.00,0.00){$\alpha$}}
\put(60.00,-5.00){\makebox(0.00,0.00){$\alpha$}}
\put(10.00,60.00){\makebox(0.00,0.00)[r]{$\theta_p$}}
\put(40.00,40.00){\makebox(0.00,0.00){$\minim_{\malesipky}$}}
\put(120.00,40.00){\makebox(0.00,0.00){$\calA_{\btb}$}}
\put(0.00,0.00){\makebox(0.00,0.00){$\minim_{\btb}$}}
\put(0.00,80.00){\makebox(0.00,0.00){$\minim_{\btb}$}}
\put(60.00,40.00){\vector(1,0){45.00}}
\put(10.00,10.00){\vector(1,1){20.00}}
\put(10.00,70.00){\vector(1,-1){20.00}}
\end{picture}}
\end{center}
\end{corollary}

\begin{remark}
{\rm\
If $\calA$ is a
non-$\Sigma$-operad, Theorem~\ref{Kukacky} holds with 
$\(\hskip -1mm -\hskip -1mm\)_{{\rm ns}}$ instead of 
$\(\hskip -1mm-\hskip -1mm \)$. 
Since $\(\hskip -1mm-\hskip -1mm\)_{{\rm ns}}$ 
has integral coefficients, 
the model $\minim_{\malesipky }$ might, as in
Example~\ref{boli_mne_hlava}, in principle exists
over the integers.  
}
\end{remark}

\begin{example}
\label{boli_mne_hlava}
{\rm\
Let us consider the free operad
$\calF(\{\mu_n,\nu_n\}_{n \geq 2}, \{p_n,q_n\}_{n \geq 1})$ generated by
\begin{itemize}
\item[(i)]
generators $\mu_n$ and $\nu_n$, $n \geq 2$,
as in Example~\ref{doufejme},
\item[(ii)]
generators
$p_n, q_n : \B^{\ot n} \to \W$ of 
degree $n-1$, $n\geq 1$, and
\item[(ii)]
generators $h_n : \B^{\otimes n} \to \W$ of
degree $n$, $n \geq 1$.
\end{itemize}

The differential of $\mu_n$ and $\nu_n$ is given by formulas~(\ref{X})
and~(\ref{Z}), the differential of $p_n$ and $q_n$ by 
formula~(\ref{Y}) (with
$f_i$ replaced by $p_i$ resp.~$q_i$, $i \geq 1$) and 
\begin{eqnarray*}
\pa (h_n) &:=& p_n -q_n +
\sum_{k=2}^n \sum_{\doubless{r_1+\cdots+r_k = n}{0 \leq s \leq k-1}} 
(-1)^{\epsilon}
\nu_k(p_{r_1} \ot \cdots \ot p_{r_s} \ot h_{r_{s+1}}  \ot q_{r_{s+2}} 
             \ot \cdots \ot q_{r_k}) +
\\
&&\hskip 1cm+ \hskip -3mm
\sum_{\doubless{i+j = n+1}{i \geq 1,j \geq 2}}
\sum_{s=0}^{n-j}
(-1)^{i+s(j+1)}
h_{i}(\id^{\ot s} \ot \mu_j \ot
\id^{\ot i-s-1}),
\end{eqnarray*}
where 
\[
\epsilon := \sum_{1 \leq i < j \leq k} r_i(r_j+1) + (r_1+\cdots
+ r_s) + k + s.
\]

An algebra over the colored operad $A(\infty)_{\malesipky}$ 
described above consists of two
$A(\infty)$-algebras ${\bf V} = (V,\partial, m_2,m_3,\ldots)$ and ${\bf W} =
(V,\partial, n_2,n_3,\ldots)$, two strongly homotopy homomorphisms
$\PP, \QQ : {\bf V} \to {\bf W}$, 
$\PP = (P_1,P_2,\ldots)$, $\QQ = (Q_1,Q_2,\ldots)$, and a
sequence $\HH = (H_1,H_2,\ldots)$ which should be interpreted as a
{\em homotopy through strongly homotopy homomorphisms\/} 
between $\PP$ and $\QQ$. This notion coincides with the homotopy
structure considered by M.~Grandis in~\cite[5.4]{grandis}.
Observe that $\pa H_1 = P_1 - Q_1$ and that, for $a,b \in V$,
\[
(\pa H_2)(a,b) =  
(-1)^a n_2 (P_1(a),H_1(b)) + n_2(H_1(a),Q_1(b))- H_1 m_2(a,b). 
\]
Writing  $m_2$ and $n_2$ as multiplications, we can translate the
above equation to
\[
(\pa H_2)(a,b) =  (-1)^a P_1(a)H_1(b) + H_1(a)Q_1(b) - H_1(ab),
\]
thus $H_1$ is `derivation homotopy up to homotopy' between $P_1$
and $Q_1$. 
}
\end{example}

\begin{remark}
{\rm\
It is well-known that $A(\infty)$-algebras are the same as
differentials on cofree coconnected coalgebras (that is, coderivations
with square zero). 
In this language, strong homotopy homomorphisms are
homomorphisms of these coalgebras commuting with the differentials. It
can be shown that the homotopies of Example~\ref{boli_mne_hlava} are
translated to ordinary (co)derivation homotopies of these maps of
coalgebras. 
}
\end{remark}

\section{Strong homotopy equivalences of algebras}
\label{Opicak_Fuk}

Let $\calA$ be an ordinary operad and $f : \B \to \W$ and $g : \W \to
\B$ be two degree zero generators. Consider the operad
\begin{equation}
\label{elmers}
\calA_{\Iso}
:=
\frac{\calA_\B * \calA_\W * \calF (f,g)}
                  {(fa_\B = a_\W \otexp fn,\ ga_\W = a_\B \otexp gn,\
fg = \id_W,\ gf = \id_\B)}
\end{equation}
describing two mutually inverse homomorphisms of $\calA$-algebras. By
our general philosophy, cofibrant resolutions of this operad will
describe homotopy equivalences of strongly homotopy $\calA$-algebras.
If $\calA$ is the operad ${\bf 1}$ for trivial algebras,
construction~(\ref{elmers}) gives the operad
\[
\Iso := {\bf 1}_\Iso = \frac{\calF(f,g)}
                  {(fg = \id_\W,\ gf = \id_\B)}
\]
describing two mutually inverse chain maps. Operad $\Iso$ was studied in
great details in~\cite{markl:ha} where we described its minimal cofibrant
resolution  as a graded colored differential operad
\begin{equation}
\label{nove_hry}
\Riso := (\calF(f_0,f_1,\ldots; g_0,g_1,\ldots), d),
\end{equation}
with two types of generators,
\begin{equation}
\label{beranek}
\def\arraystretch{1.25}
\begin{array}{rl}
\mbox{(i)}&\mbox{% 
\hskip -2mm 
generators $\seq{f_k}{k \geq 0}$, $\deg(f_k)=k$,} 
\left\{
\begin{array}{l}
\mbox{$f_k : \B \to \W$ if $k$ is even,} 
\\
\mbox{$f_k : \B \to \B$ if $k$ is odd,}
\end{array}
\right.
\\
\mbox{(ii)}&\mbox{% 
\hskip -2mm
generators
$\seq{g_k}{k \geq 0}$, $\deg(g_k)=k$,}
\left\{
\begin{array}{l} 
\mbox{$g_k : \W \to \B$ if $k$ is
even,}
\\
\mbox{$g_k : \W \to \W$ if $k$ is odd.}
\end{array}
\right.
\end{array}
\end{equation}
The differential $d$ was given by
\begin{equation}
\label{virzinko}
\def\arraystretch{1.4}
\begin{array}{ll}
d f_0\widedef 0, 
& 
d g_0 \widedef 0, 
\\
d f_1 \widedef g_0f_0 - 1,   
& 
d g_1 \widedef f_0g_0 - 1    
\end{array}
\end{equation}
and, on remaining generators, by the formula
\begin{eqnarray}
\nonumber
df_{2m}   &:=& \sum_{0 \leq i < m}(f_{2i}f_{2(m-i)-1}-
g_{2(m-i)-1}f_{2i}),\ m \geq 0,
\\
\label{pejsek_a_kocicka}
df_{2m+1} &:=& \sum_{0 \leq j \leq m} g_{2j}f_{2(m-j)} -
               \sum_{0 \leq j < m} f_{2j+1}f_{2(m-j)-1},\ m \geq 1,
\\
\nonumber 
dg_{2m}   &:=& \sum_{0 \leq i < m}(g_{2i}g_{2(m-i)-1}-
f_{2(m-i)-1}g_{2i}),\ m \geq 0,
\\
\nonumber 
dg_{2m+1} &:=& \sum_{0 \leq j \leq m} f_{2j}g_{2(m-j)} -
               \sum_{0 \leq j < m} g_{2j+1}g_{2(m-j)-1},\ m \geq 1.
\end{eqnarray}
The above formulas can be written in a shorter form by
introducing `formal generators'
\begin{equation}
\label{bojim_se_porad}
\def\arraystretch{1.5}
\begin{array}{ll}
\fb := f_0 + f_2 + f_4 + \cdots : \B \to \W, &
\hb := f_1 + f_3 + f_5 + \cdots : \B \to \B,
\\
\gb := g_0 + g_2 + g_4 + \cdots : \W \to \B, &
\lb := g_1 + g_3 + g_5 + \cdots : \W \to \W.
\end{array}
\end{equation} 
Then $\Riso = \calF(\fb, \gb, \hb, \lb)$ with 
the differential
\[
d\fb = \fb\hb - \lb\fb,\   d\hb = \gb\fb - \hb\hb -1_\B,\
d\gb = \gb\lb - \hb\gb \mbox { and }   d\lb = \fb\gb - \lb\lb -1_\W. 
\]
\noindent 
{\boldface Warning.}
We will use the above abbreviation very often, but we
shall always keep in mind that each formula of this type in fact represents
infinitely many formulas for homogeneous parts. Observe also
that in this section the symbols $f_0,f_1,\ldots$ have 
different meaning than in Example~\ref{doufejme} and related examples
and applications.

We will need the following definition.

\begin{definition}
\label{snezenky}
Let $m \geq 1$. A polarization is a choice of $\Sigma_m$-equivariant
morphisms 
$\(\fb\) \in \Riso(\otexp{\B}m, \otexp{\W}m)$,
$\(\gb\) \in \Riso(\otexp{\W}m, \otexp{\B}m)$,
$\(\hb\) \in \Riso(\otexp{\B}m, \otexp{\B}m)$ and
$\(\lb\) \in \Riso(\otexp{\W}m, \otexp{\W}m)$ such that
\begin{itemize}
\item[(i)]
$\(\fb\)_0 = \otexp{f_0}m$ and $\(\gb\)_0 = \otexp{g_0}m$ (the
subscript $0$ denotes the degree zero part), and
\item[(ii)] 
the following differential equations in $\Riso$ are satisfied:
\[
\def\arraystretch{1.4}
\begin{array}{ll}
d\(\fb\) = \(\fb\)\(\hb\) - \(\lb\)\(\fb\),
& 
d\(\hb\) = \(\gb\)\(\fb\) - \(\hb\)\(\hb\) -1_\B,
\\
d\(\gb\) = \(\gb\)\(\lb\) - \(\hb\)\(\gb\),
& 
d\(\lb\) = \(\fb\)\(\gb\) - \(\lb\)\(\lb\) -1_\W.
\end{array}
\]
\end{itemize}
\end{definition}

The existence of a polarization easily follows from the acyclicity of
the resolution $\Riso$ proved in~\cite{markl:ip}. 
See also Remark~\ref{polari}.

Let $X$ be a collection. In the main theorem of this section we 
consider the free operad 
\[
(\calF(X_\B; X_\W; \{f_k\}_{k \geq 0}; \{g_k\}_{k \geq 0},
\{X^{f_k}\}_{k \geq 0}; \{X^{g_k}\}_{k \geq 0}),
\]
 generated by
\begin{itemize}
\item[(i)]
two copies $X_\B$ and $X_\W$ of $X$ interpreted as in
Theorem~\ref{Dan1}, 
\item[(ii)] 
the generators $f_k$ and $g_k$, $k \geq 0$, as in~(\ref{beranek}) and
\item[(iii)]
two copies $X^{f_k}$ and $X^{g_k}$ of the suspension $\suspit {k+1}
X$, $k \geq 0$.
\end{itemize}
A generator $x^{f_k} \in X^{f_k}$ is, for $x \in X(n)$, interpreted as
a morphism
$\otexp {\B}n \to \W$ if $k$ is even and as a morphism $\otexp {\B}n \to \B$
if $k$ is odd. Similarly, $x^{g_k} \in X^{g_k}$ is interpreted as a morphism
$\otexp {\W}n \to \B$ if $k$ is even and as a morphism $\otexp {\W}n \to \W$
if $k$ is odd.

\begin{theorem}
\label{asi_do_te_Ameriky_neodjedu}
Let $(\calA,d_\calA)$ be an (ordinary) operad satisfying~(\ref{Jitka})
and let $\rho :\minim \to (\calA,d_\calA)$, 
$\minim = (\calF(X),d_{\minim})$, 
be its minimal model.
Then the minimal model $\minim_{\Iso}$ of the colored 
operad $\calA_\Iso$ is of the form
\[
\beta:
(\calF(X_\B; X_\W; \{f_k\}_{k \geq 0}; \{g_k\}_{k \geq 0},
\{X^{f_k}\}_{k \geq 0}; \{X^{g_k}\}_{k \geq 0}), D) \longrightarrow 
(\calA_{\Iso},d)
\]
with $\beta$ given by
$\beta(x_{\B}) = j_{\B}(\rho(x))$, $\beta(x_{\W}) = j_{\W}(\rho(x))$
for $x \in X$, 
$\beta(f_0) = f$, $\beta(g_0) = g$, while $\beta$ is trivial on the
remaining generators.

The differential $D$ is on $f_k$ and $g_k$, $k \geq
0$, given by formulas~(\ref{virzinko}) and~(\ref{pejsek_a_kocicka}).
For $n \geq 2$ and $x \in X(n)$,  $D(x_{\B}) = j_{\B}(d_{\minim}(x))$, 
$D(x_{\W}) = j_{\W}(d_{\minim}(x))$ and, in the shorthand 
of~(\ref{bojim_se_porad}), 
\begin{eqnarray*}
D(x^\fb) &=& \fb x_\B - x_\W \(\fb\) + \fb x^\hb -(-1)^{x} x^\fb \(\hb\) 
            - \lb x^\fb -x^\lb \(\fb\) +  \omega_\fb,
\\
D(x^\gb) &=& \gb x_\W - x_\B \(\gb\) + \fb x^\lb -(-1)^{x}  x^\gb \(\lb\) 
            - \hb x^\gb -x^\hb \(\gb\) + \omega_\gb,
\\
D(x^\hb) &=&  (-1)^x x_\B \(\hb\) - \hb x_\B- \hb x^\hb  +(-1)^x
            x^\hb \(\hb\) - \gb x^{\fb} + x^\gb \(\fb\) + \omega_\hb,
\\
D(x^\lb) &=&  (-1)^x x_\W \(\lb\) - \lb x_\W - \lb x^\lb +(-1)^x x^\lb \(\lb\) 
                - \fb x^{\gb} + x^\fb \(\gb\) + \omega_\lb. 
\end{eqnarray*}
where $\omega_\fb$, $\omega_\gb$, $\omega_\hb$ and $\omega_\lb$ depend
linearly on $x$ and are elements of the ideal generated by
$X^{f_k}(\<n)$, $X^{g_k}(\<n)$, $k \geq 0$.
\end{theorem}

\noindent 
{\boldface Proof.}
Assume that the differential $D$ exists. The map $\beta$ then decomposes as
\[
\minim_\Iso
\stackrel{\beta_1}{\lra}
\left(\frac{
    \calF(X_\B;X_\W; \fb, \gb, \hb, \lb)}I,d
\right)
\stackrel{\beta_2}{\lra} \calA_\Iso,
\]
where the ideal $I$ is generated by
\[
\fb x_\B - x_\W \(\fb\),\ \gb x_\W - x_\B \(\gb\),\
 (-1)^x x_\B \(\hb\) - \hb x_\B, \mbox { and } 
 (-1)^x x_\W \(\lb\) - \lb x_\W,\ x\in X.
\]
The definitions of $\beta_1$ and $\beta_2$ is obvious. 
To prove that $\beta_1$ is a quism, consider the filtration of
$\minim_\Iso$ given by the number of generators from $X^{f_k}$,
$X^{g_k}$, $k \geq 0$ and use a standard spectral sequence
argument. The nature of $\beta_2$ is similar to that of the map
$\beta_\btb$ in Corollary~\ref{snad_to} and the same 
arguments as in its proof show that $\beta_2$ is a quism. The tails can be
constructed by a general nonsense. We leave the details to the reader.%
\qed

\begin{remark}
\label{polari}
{\rm\
There is a very explicit, though rather involved, formula for
the polarization of Definition~\ref{snezenky}. 
If we do not demand the equivariance,
there exists (in fact many)
a polarization with {\em integral\/} coefficients. An example 
is given, for $m=2$, by
\[
\def\arraystretch{1.4}
\begin{array}{ll}
\(\fb\) := \sum_{i \geq 0} \fb \otexp{(\gb \fb)}i \ot f_{2i},
& 
\(\gb\) := \sum_{i \geq 0} \gb \otexp{(\fb \gb)}i \ot g_{2i},
\\
\(\hb\) := \hb \ot \id + 
             \sum_{i \geq 1} \otexp{(\gb \fb)}i \ot h_{2i-1},
& 
\(\lb\) := \lb \ot \id + 
             \sum_{i \geq 1} \otexp{(\fb \gb)}i \ot l_{2i-1}.
\end{array}
\]
and by similar formulas for $m > 2$. These formulas may be used for
integral minimal models of $\calA_\Iso$ for non-$\Sigma$ operads.
}
\end{remark}

We believe that, for $\calA$ quadratic Koszul, there exists a closed
formula for the `tails' of the differential $D$. It might be a
challenging, though very involved, exercise to write this formula 
when $\calA =\Ass$, the operad describing associative algebras.

\section{Final remarks and challenges}
\label{Misa}

In this section, ``sh'' will abbreviate ``strongly (or strong) homotopy.''

\noindent 
{\boldface The `category' of sh{} algebras.}
Let $\calA$ be an ordinary operad. Recall that a {\em sh{}
$\calA$-algebra\/} is an algebra over the minimal model 
$\minim = (\calF(X),d_{\minim})$ of $\calA$. {\em A sh{} homomorphism of sh{}
$\calA$-algebras\/} is then a colored algebra over the 
minimal model $\minim_\btb$
of $\calA_\btb$. See Example~\ref{doufejme} for the case $\calA =
\Ass$. 

Let $\PP : \minim_{\B \to \W} \to \End_{U,V}$ and $\QQ : \minim_{\W
\to \T} \to \End_{V,W}$, where $\T$ is a \underline{{\tt T}}hird color, be
two sh{} homomorphisms. Their `composition' $\QQ \circ \PP :
\minim_{\B \to \T} \to \End_{U,W}$ is given by an appropriate map
\begin{equation}
\label{Buxheimer_MS}
\Xi : \minim_{\B \to \T} \to \minim_{\B \to \W} * \minim_{\W \to \T}
\end{equation}
($*$ denotes the free product), whose existence is guaranteed 
by a general nonsense.

A topological counterpart of this situation was studied by Boardman and
Vogt. It turned out that the composition above  need not
be associative, thus topological sh{}-algebras do not
form a honest category, but only a hazier object called weak
Kan category~\cite[Theorem~4.9]{boardman-vogt:73}.

We are sure that the situation in algebra is much better, because
minimal models are `reasonably' functorial. 
So we believe that the following problem has an affirmative answer.

\begin{problem}
The map $\Xi$ in~(\ref{Buxheimer_MS}) can be constructed in such a way that
the composition $\circ$ is associative, that is, sh{} algebras and
their sh{} homomorphisms form a honest category.
\end{problem}

Since sh{} algebras over quadratic Koszul operads may be interpreted
as certain dg cofree coconnected coalgebras and their maps as
homomorphisms of these coalgebras, for these operads the map $\Xi$
inducing an associative composition clearly exists.

\vskip 2mm
\noindent 
{\boldface Ideal perturbation lemma.}
In~\cite{markl:ip} we proved that a
proper formulation of the homological perturbation lemma needs a 
good notion of
chain homotopy equivalences. These `good homotopy equivalences' were
defined as representations of the colored operad $\Riso$ of~(\ref{nove_hry})
and called {\em strong\/} homotopy equivalences.
Our theory was, however, tailored for chain complexes 
with no additional algebraic structure.

Constructions of Section~\ref{Opicak_Fuk} provide a similar notion 
for chain complexes with an additional algebraic structure:

\begin{definition}
\label{zase_spatna_disketa}
Let $A:\minim \to \End_{U}$ and $B : \minim \to \End_V$ be two sh{}
$\calA$-algebras. 
A strong homotopy equivalence between $A$ and $B$ is
given by an operad action $S : \minim_\Iso \to \End_{U,V}$ such that
$S\circ j_\B = A$ and $S \circ j_\W = B$, where
\[
j_\B : \minim \hookrightarrow \minim_\Iso \hookleftarrow \minim : j_\W
\]
are obvious inclusions.
\end{definition}

\begin{problem}
Formulate an `ideal perturbation lemma' for
chain complexes with an additional algebraic structure, using sh{}
equivalences of Definition~\ref{zase_spatna_disketa}.
\end{problem}

\noindent
{\boldface Homotopies through homomorphisms.}
While it is quite obvious what a ``homotopy through homomorphisms'' 
for maps of topological monoids should be, the situation in algebra is
much subtler.

Suppose we have two sh{} $\calA$-algebras $A: \minim \to \End_U$, $B:
\minim \to \End_V$ and two sh{} homomorphisms
$\PP,\QQ : A \to B$ given by actions $\PP,\QQ : \minim_\btb \to
\End_{U,V}$. 
These data clearly define a representation of the outer square
of the diamond
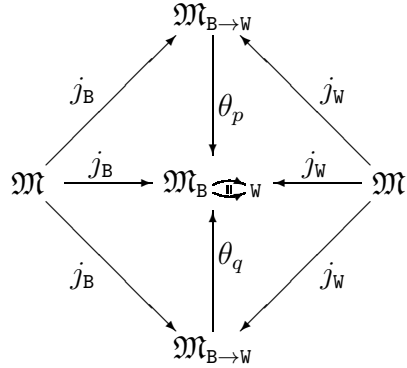
\begin{center}

{% Picture saved by xtexcad 2.4
\unitlength=.7pt
\begin{picture}(180.00,185.00)(0.00,0.00)
\put(-10,0){
\put(30.00,40.00){\makebox(0.00,0.00){$j_{\B}$}}
\put(40.00,100.00){\makebox(0.00,0.00){$j_{\B}$}}
\put(30.00,140.00){\makebox(0.00,0.00){$j_{\B}$}}
\put(10.00,80.00){\vector(1,-1){70.00}}
\put(10.00,100.00){\vector(1,1){70.00}}
\put(0.00,90.00){\makebox(0.00,0.00){$\minim$}}
}

\put(5,0){
\put(140.00,100.00){\makebox(0.00,0.00){$j_{\W}$}}
\put(150.00,40.00){\makebox(0.00,0.00){$j_{\W}$}}
\put(150.00,140.00){\makebox(0.00,0.00){$j_{\W}$}}
\put(170.00,100.00){\vector(-1,1){70.00}}
\put(170.00,80.00){\vector(-1,-1){70.00}}
\put(180.00,90.00){\makebox(0.00,0.00){$\minim$}}
}

\put(10.00,90.00){\vector(1,0){47.00}}
\put(100.00,50.00){\makebox(0.00,0.00){$\theta_q$}}
\put(100.00,130.00){\makebox(0.00,0.00){$\theta_p$}}
\put(90.00,90.00){\makebox(0.00,0.00){$\minim_{\malesipky}$}}
\put(90.00,0.00){\makebox(0.00,0.00){$\minim_{\btb}$}}
\put(90.00,180.00){\makebox(0.00,0.00){$\minim_{\btb}$}}
\put(90.00,170.00){\vector(0,-1){65.00}}
\put(90.00,10.00){\vector(0,1){65.00}}
\put(170.00,90.00){\vector(-1,0){47.00}}

\end{picture}}
\end{center}

\begin{definition}
\label{Scastjem}
A homotopy through sh{} homomorphisms 
between $\PP$ and $\QQ$ is
an operad action $\minim_\malesipky \to \End_{U,V}$
extending the data above. We write $\PP\sim\QQ$
if there exists a sh{} homotopy between $\PP$ and $\QQ$.
\end{definition}

As we saw in Example~\ref{boli_mne_hlava}, this notion generalizes the
homotopy structure on the category of $A(\infty)$-algebras considered
in~\cite[5.4]{grandis}.  
The relation $\sim$ is indeed an equivalence:

\begin{proposition}
The relation $\sim$ is an equivalence on the set of sh{} homomorphisms
of sh{} $\calA$-algebras.
\end{proposition}

\noindent 
{\boldface Proof.}
The reflexivity ($\PP \sim \PP$) and symmetry ($\PP \sim
\QQ \Longrightarrow \PP \sim
\QQ$) is obvious. If $\PP \sim \QQ$ and $\QQ \sim \TT$ 
then the homotopy between $\PP$ and $\TT$ can be
constructed by standard homological methods (general nonsense). We
leave the details to the reader.%
\qed

\noindent 
{\boldface Possible generalizations.}
The main results of the paper (Theorems~\ref{Dan1},~\ref{Kukacky} 
and~\ref{asi_do_te_Ameriky_neodjedu}) 
seem to admit a generalization to resolutions of colored
operads $\calA_\frakD$ describing $\frakD$-diagrams of $\calA$-algebras for
arbitrary $\frakD$. 

We start with  a cofibrant resolution ${\cal R}_\frakD = (\calF(F),d_\frakD)$
of the diagram $\frakD$ considered as a
$V$-colored operad, where $V$ is the set of vertices of $\frakD$. 
The second ingredient is the minimal model 
$\minim = (\calF(X),d_{\minim})$ of $\calA$. We believe that it is
possible to prove:

\begin{conjecture}
There exists a cofibrant model $\minim_\frakD$ 
of the colored operad $\calA_{\frakD}$ of the form
\[
\minim_\frakD =
(\calF(X\times V \sqcup F \sqcup X^F),D),
\]
where $X \times \{v\} \subset X \times V$ is, for $v \in V$, a copy of $X$
``concentrated'' in the color $v$ and elements of $X^f$ are considered as
multilinear maps 
with the same source and the same target as $f \in F$. 
The differential $D$ on
$X\times V$ and on $F$ is given by the inclusions
\[
j_\frakD : {\cal R}_\frakD \hookrightarrow \minim_\frakD \mbox { and }
j_v : \minim \hookrightarrow \minim_\frakD,\
v \in V,
\] 
while the differential of $x^f$, for $f \in F$ and $x \in X(n)$, $n
\geq 2$, decomposes as
\[
D(x^f) = {\it Pr} + \omega_x.
\]
The principal part ${\it Pr}$ further
decomposes as ${\it Pr} = {\it Pr}_1 + {\it Pr}_2$, where ${\it Pr}_1$
is `ideologically' the commutator of $x$ and $f$ and ${\it Pr}_2$ is given
by a polarization of $d_\frakD(f)$. The tail $\omega_x$ belongs to the
ideal generated by $X^f(\<n)$, $f \in F$.
The model $\minim_\frakD$ is minimal if and only if ${\cal R}_\frakD$
is. 
\end{conjecture}

We leave as an exercise to check that all models described in this
paper are of the above form (in all cases $V = \{\B,\W\}$).

%\bibliography{b}

\begin{thebibliography}{1}

\bibitem{boardman-vogt:73}
J.M. Boardman and R.M. Vogt.
\newblock {\em Homotopy Invariant Algebraic Structures on Topological Spaces}.
\newblock Springer-Verlag, 1973.

\bibitem{ginzburg-kapranov:DMJ94}
V.~Ginzburg and M.M. Kapranov.
\newblock Koszul duality for operads.
\newblock {\em Duke Math. Journal}, 76(1):203--272, 1994.

\bibitem{grandis}
M.~Grandis.
\newblock On the homotopy structure of strongly homotopy associative algebras.
\newblock {\em Jour. Pure Appl. Algebra}, 134(1):15--81, 1999.

\bibitem{lada-stasheff:IJTP93}
T.~Lada and J.D. Stasheff.
\newblock Introduction to sh {L}ie algebras for physicists.
\newblock {\em International Journal of Theoretical Physics}, 32(7):1087--1103,
  1993.

\bibitem{markl:JPAA92}
M.~Markl.
\newblock A cohomology theory for {$A(m)$-algebras} and applications.
\newblock {\em Jour. Pure Appl. Algebra}, 83:141--175, 1992.

\bibitem{markl:zebrulka}
M.~Markl.
\newblock Models for operads.
\newblock {\em Communications in Algebra}, 24(4):1471--1500, 1996.

\bibitem{markl:ha}
M.~Markl.
\newblock Homotopy algebras are homotopy algebras.
\newblock Preprint {\tt math.AT/9907138}, July 1999.

\bibitem{markl:ip}
M.~Markl.
\newblock Ideal perturbation lemma.
\newblock Preprint {\tt math.AT/0002130}, to appear in Communications in
  Algebra.

\bibitem{may:1972}
J.P. May.
\newblock {\em The Geometry of Iterated Loop Spaces}, volume 271 of {\em
  Lecture Notes in Mathematics}.
\newblock Springer-Verlag, 1972.

\end{thebibliography}

\vskip3mm
\catcode`\@=11
\noindent
Mathematical Institute of the Academy, \v Zitn\'a 25, 115 67
Praha 1, The Czech Republic,\hfill\break\noindent
email: {\tt
markl@math.cas.cz}\hfill\break\noindent

\vfill

\hfill {\tt \jobname.tex}

\end{document}